\documentclass{amsart}

\newtheorem{theorem}{Theorem}[section]

\theoremstyle{definition}

\theoremstyle{remark}

\numberwithin{equation}{section}

\begin{document}

\title{Dynamics of Electrostatic MEMS Actuators}
\newcommand{\bfR}{{\Bbb R}}\newcommand{\ep}{\epsilon}
\newcommand{\bfC}{{\Bbb C}}
\newcommand{\bfH}{{\Bbb H}}\newcommand{\coker}{\mbox{coker}}\newcommand{\image}{\mbox{image}}
\newcommand{\bfZ}{{\Bbb Z}}\newcommand{\bfP}{{\Bbb P}}
\newcommand{\bA}{{{\bf A}}}\newcommand{\bF}{{\bf F}}
\newcommand{\De}{{\Delta}}
\newcommand{\sech}{\mbox{ sech}}
\newcommand{\DD}{\mbox{ D}}
\newcommand{\F}{\vec{F}}
\newcommand{\Ph}{\vec{\phi}}\newcommand{\x}{{\bf x}}\newcommand{\y}{{\bf y}}\newcommand{\f}{{\bf f}}\newcommand{\z}{{\bf z}}
\newcommand{\A}{\vec{A}}
\newcommand{\bfph}{{\bf \phi}}
\newcommand{\lm}{\lambda}

\newcommand{\dd}{\mbox{d}}\newcommand{\D}{\cal O}
\newcommand{\ee}{\end{equation}}
\newcommand{\be}{\begin{equation}}
\newcommand{\bea}{\begin{eqnarray}}
\newcommand{\eea}{\end{eqnarray}}
\newcommand{\ii}{\mbox{i}}\newcommand{\e}{\mbox{e}}\newcommand{\rr}{{\mbox{r}}}
\newcommand{\pa}{\partial}\newcommand{\Om}{\Omega}\newcommand{\om}{\omega}
\newcommand{\vep}{\varepsilon}
\newcommand{\di}{\displaystyle}
\newcommand{\nn}{\nonumber}
\newcommand{\V}{V_{\mbox{Dpi}}}\newcommand{\Vc}{V_{\mbox{c}}}\newcommand{\VV}{V_{{\rm{\mbox{Dpi}}}}}
\newcommand{\xd}{x_{\mbox{Dpi}}}\newcommand{\xxd}{x_{\rm{\mbox{Dpi}}}}
\newcommand{\xs}{x_{\mbox{s}}}\newcommand{\xxs}{x_{\rm{\mbox{s}}}}\newcommand{\ts}{t_{\mbox{s}}}\newcommand{\tts}{t_{\rm{\mbox{s}}}}
\newcommand{\tp}{t_{\mbox{p}}}\newcommand{\ttp}{t_{\rm{\mbox{p}}}}\newcommand{\tc}{t_{\mbox{c}}}\newcommand{\ttc}{t_{\rm{\mbox{c}}}}

\author{Yisong Yang}
\address{Department of Mathematics, Polytechnic Institute of New York University, Brooklyn, New York 11201, USA}
\email{yyang@math.poly.edu}

\author{Ruifeng Zhang}
\address{Institute of Contemporary Mathematics and School of Mathematics, Henan University, Kaifeng, Henan 475004, PR China}
\email{zrf615@henu.edu.cn}
\author{Le Zhao}
\address{School of Mathematics, Henan University, Kaifeng, Henan 475004, PR China}





\begin{abstract} Electrostatic actuators are simple but important switching devices for MEMS applications. Due to the difficulties associated with
the electrostatic nonlinearity, precise mathematical description is often hard to obtain for the dynamics of these actuators.
Here we present two sharp theorems concerning the dynamics of an undamped electrostatic actuator with one-degree of freedom, subject to linear and nonlinear elastic forces, respectively.
We prove that both situations are characterized by the onset of one-stagnation-point periodic response below a well-defined pull-in voltage and a finite-time touch-down or collapse of the actuator above this pull-in voltage.
In the linear-force situation, the stagnation level, pull-in voltage, and pull-in coordinate of the movable electrode may all be determined explicitly,
following the recent work of Leus and Elata based on numerics. Furthermore, in the nonlinear-force situation, the stagnation level, pull-in voltage, and pull-in coordinate may be described completely in terms of the electrostatic and mechanical parameters
of the model so that they approach those in the linear-force situation monotonically in the zero nonlinear-force limit.
\end{abstract}

\maketitle



\section{Introduction}

It is well known that electrostatic actuators play a basic and important role in microelectromechanical systems (MEMS) technology as one of the most commonly used microcircuit components
 and mathematical modeling and numerical simulation are essential for understanding
various qualitative and quantitative properties of such devices \cite{LY,PB}. However, due to the inverse-square type nonlinearity
in the governing equations associated with the Coulomb electrostatic forces, it is a challenging task for analysts to gather a precise picture of the 
mathematical behavior of the systems and numerical simulations have been resorted to in the investigations of the modeling of actuators. So far, overwhelmingly, mathematical studies have been
concentrated on equations governing steady state solutions \cite{DFG,EGG,GW,LY,M,PB} or their parabolic extensions \cite{EGG,YZ} and little work has been done concerning the dynamical
behavior of the problem which is of a wave-equation
nature and of obvious importance. (In \cite{CKL,Guo}, some existence, uniqueness, and continuous dependence results are obtained for
the evolution problems related to the MEMS wave equations. We thank the referee for informing us these developments.) In this paper, as a first attempt, we shall study the dynamics of an undamped actuator of one-degree of freedom described by a second-order ordinary differential equation, as a simplified
model arising from the original spatially dependent wave equation description \cite{PB}.
The simplicity of the model enables us to obtain
 two sharp theorems concerning the dynamics of the actuator, subject to linear and nonlinear elastic forces, respectively,
which roughly say that in both situations the dynamic response of the actuator
under an applied voltage $V$ is characterized by the onset of one-stagnation-point periodic motion when $V$ is below a well-defined pull-in voltage and a finite-time touch-down or collapse of the actuator happens when $V$ is above this pull-in voltage. More precisely,
in the linear-force situation, the stagnation level, pull-in voltage, and pull-in coordinate of the movable electrode may all be determined explicitly,
following the recent work of Leus and Elata \cite{LE} based on numerics, while, in the nonlinear-force situation, the stagnation level, pull-in voltage, and pull-in coordinate may be described completely in terms of the electrostatic and mechanical parameters
of the model so that these quantities approach those in the linear-force situation monotonically in the zero nonlinear-force limit.
It is interesting that, under a fairly natural convexity assumption, our method for the linear-force situation is completely applicable to the nonlinear situation so that an understanding of the global dynamics of the undamped
actuator is achieved.

 The rest of the paper is organized as follows. In Section 2, we review the physical description and equation of motion of the actuator subject to a linear elastic force
and the work of Leus and Elata \cite{LE}, which motivates our study to follow. In Section 3, we use the energy conservation law as a first integral to establish the existence of
a stagnation time and stagnation level for the motion when the applied voltage is below a critical value. In Section 4, we use the results of Section 3 to show that the motion
with a stagnation point actually gives rise to a symmetric periodic motion of the actuator.
In Section 5, we show that when the applied voltage is above the critical value identified in \cite{LE} and Section 4, the solution is non-periodic and leads to a touch-down or collapse 
situation of the device. In Section 6, we summarize our study as a theorem which describes the full dynamical behavior of the actuator and confirms the results stated in \cite{LE},
when the elastic force is linear. In Section 7, we extend our method to study a more general model when the elastic force is cubicly nonlinear which is widely useful in applications
(cf. \cite{CO,FM} and references therein). We conclude by indicating how our method may be used to tackle models under more general nonlinear elastic forces.
In Section 8, we spell out some simple but generic conditions under which the dynamic touch-down or collapse phenomenon occurs universally
 as a result of a sufficiently high applied voltage.
\medskip

\section{ Physical model and governing equation} 

Following \cite{LE}, the parallel-plate actuator of one-degree of freedom
under study consists of a top electrode of mass $m$ and area $A$ that is suspended on a 
linear elastic spring with the Hooke constant $k$, above a fixed bottom electrode that is coated with a dielectric layer of thickness $d_0\geq0$. The initial gap 
or distance between the top electrode
and the dielectric is $g>0$. The fixed bottom electrode is electrically grounded and a step-function type voltage will be applied to the top electrode so that it is zero before time $t=0$
and assumes the constant value $V>0$ when $t>0$. Such an application results in the one-dimensional motion of the top electrode away from its initial equilibrium position 
with the coordinate $x=0$ to that with $x>0$ so that the Hamiltonian of the motion, which is the sum of
the kinetic energy, elastic potential arising from a linear elastic force, and electrostatic potential of the deformable capacitor subject to the voltage $V$, is given by the expression \cite{LE}
\be \label{H1}
H=\frac12 m{\dot{x}}^2+\frac12 k x^2-\frac{\varepsilon_0 AV^2}{2\left(\frac{d_0}{\varepsilon_\rr}+g -x\right)},
\ee
where $\dot{x}=\dd x/\dd t$ is the velocity of the moving top plate, $\vep_0$ the permittivity of the free space, and $\vep_\rr$ the relative permittivity of the dielectric.

In order to simply the notation, we may conveniently adopt the normalized variables and parameters \cite{LE} as follows
\be \label{V}
 \frac xg\mapsto x,\quad \frac H{kg^2}\mapsto H,\quad\sqrt{\frac km}t\mapsto t,\quad \frac{\vep_0 A V^2}{kg^3}\mapsto V^2,\quad\frac{d_0}{g\vep_\rr}=\xi.
\ee
Hence the Hamiltonian (\ref{H1}) takes the form 
\be\label{H2}
H=\frac12{\dot{x}}^2+\frac12 x^2-\frac12\frac{V^2}{(\xi+1-x)}.
\ee
As a consequence, the momentum equation governing the motion of the top plate-shaped electrode is \cite{LE}
\be\label{Eq}
{\ddot{x}}=-\frac{\pa H}{\pa x}=-x+\frac12\frac{V^2}{(\xi+1-x)^2},\quad t>0.
\ee

Since the actuator is free before the onset of the positive voltage $V>0$ starting at $t=0$, the initial condition for the coordinate variable of the top electrode
is given by
\be\label{IC}
x(0)=0,\quad {\dot{x}}(0)=0.
\ee

The simple initial value problem consisting of (\ref{Eq}) and (\ref{IC}) and describing the dynamic response of an undamped actuator \cite{LE}, subject to a step-function type applied
voltage and a linear elastic force, will be our main focus in 
the first part of our study here. Although the nonlinearity arising from the Coulomb electrostatic force appears to be singular at $x=1+\xi$, the touch-down 
or collapse actually occurs
at
\be 
x=1,
\ee
due to the fact that the bottom electrode is coated with a dielectric layer of thickness $d_0\geq0$ and relative permittivity $\vep_\rr$, which results in an effective insulator layer of
the `thickness' $\xi\geq0$ given in (\ref{V}) in the normalized units.

Based on numerical simulations, Leus and Elata \cite{LE} find that there is a critical voltage, called the `dynamic pull-in voltage', $\V$, so that 
below $\V$ the dynamic response of the actuator is periodic, above $\V$, the response is non-periodic so that the motion is characterized by a positive velocity which results
in a finite-time touch-down or collapse, and at the critical voltage $\V$ the movable electrode converges to an unstable equilibrium 
position as $t\to\infty$, called the `dynamic pull-in' position,
$\xd>0$. In order to calculate the dynamic pull-in voltage and dynamic pull-in position, Leus and Elata \cite{LE} investigated the periodic solutions which, according to their numerical simulations, would
oscillate between the initial state $x=0$ and the state of the maximum displacement, $\xs$, where the velocity vanishes, called the stagnation displacement. In order to find $\xs$, they use
the energy conservation law
\be 
H|_{t>0}=H |_{t=0}=-\frac12\frac{V^2}{\xi+1}
\ee
to get the first integral
\be \label{FI}
{\dot{x}}^2+ x^2 -\frac{xV^2}{(\xi+1-x)(\xi+1)}=0
\ee
for any solution at $t>0$. Hence, the stagnation position $\xs$ may be obtained by setting ${\dot{x}}=0$ in (\ref{FI}) to yield
\be 
(\xi+1)\xs^2-(\xi+1)^2\xs+V^2=0.
\ee
These equation has two positive roots, one is below $(\xi+1)/2$ and the other above $(\xi+1)/2$; both are below $\xi+1$. Leus and Elata dismiss the larger root as non-physical and choose
 the smaller root to determine the stagnation position,
\be \label{xs}
\xs=\frac{(\xi+1)^2-\sqrt{(\xi+1)^4-4(\xi+1)V^2}}{2(\xi+1)}.
\ee
The formula (\ref{xs}) clearly indicates that the stagnation position increases with regard to the applied voltage until it fails to exist when the argument in the square root of (\ref{xs}) becomes negative.
As a consequence, they conclude that the dynamic pull-in voltage and dynamic pull-in position are given by the explicit formulas
\be 
\V=\frac12(1+\xi)^{\frac32},\quad \xd=\frac12(1+\xi).
\ee

In the first part of this paper consisting of Sections 3--6, we aim at establishing the above global dynamic picture for an undamped parallel-plate actuator
subject to a linear elastic force rigorously. In the second part of the paper (Section 7), we will use the method developed in detail in the previous sections to obtain the same dynamic picture
when the actuator is under a cubicly nonlinear elastic force. In doing so, we will see that our method may readily be extended to investigate more general models.

\section{Use of the first integral}

To proceed, we rewrite the conservation law (\ref{FI}) as
\be\label{2.1}
{\dot{x}}^2=\frac{V^2}{(\xi+1)}\frac{x}{(\xi+1-x)}-x^2.
\ee
We see from (\ref{2.1}) that any solution $x(t)$ of the initial value problem consisting of (\ref{Eq}) and (\ref{IC}) must remain nonnegative-valued for all $t$,
\be\label{1.5}
x(t)\geq 0.
\ee

For convenience, we denote the right-hand side of (\ref{2.1}) by $f(x)$. Then it takes the factored form
\be\label{2.2}
f(x)=\frac{x}{(\xi+1-x)}\bigg(x^2-(\xi+1) x +\frac{V^2}{(\xi+1)}\bigg)\equiv \frac{x}{(\xi+1-x)}q(x).
\ee
We see that, when the voltage $V$ is low such that
\be\label{2.3}
V^2<\frac14 (\xi+1)^3,
\ee
the quadratic polynomial $q(x)$ has two real roots, say $x_1$ and $x_2$, given by the expressions
\be
x_{1,2}=\frac12(\xi+1)\mp\frac12\sqrt{(\xi+1)^2-\frac{4V^2}{(\xi+1)}}.\label{2.4}\\
\ee
It is important to observe the property 
\be\label{2.5}
0<x_1<x_2<\xi+1,
\ee
which allows us to express $f(x)$ as
\be\label{2.6}
f(x)=\frac{x}{(\xi+1-x)}(x_1-x)(x_2-x).
\ee
As a consequence, we see that the only physically acceptable solution has to satisfy the boundedness condition
\be\label{2.7}
0\leq x(t)\leq x_1,
\ee
otherwise it would violate the non-negativity of the right-hand side of (\ref{2.1}).

With (\ref{2.7}) in mind, we decompose (\ref{2.1}) into two pieces,
\be\label{2.8}
\dot{x}=\sqrt{\frac{x}{(\xi+1-x)}(x_1-x)(x_2-x)},\quad t>0,
\ee
and
\be\label{2.9}
\dot{x}=-\sqrt{\frac{x}{(\xi+1-x)}(x_1-x)(x_2-x)},\quad t>0,
\ee
which may be viewed as two first integrals of (\ref{Eq}). It is a comfort to note that the initial condition $x(0)=0$ automatically implies the desired
initial condition for velocity, $\dot{x}(0)=0$, stated in (\ref{IC}).

Using the condition $x(t)\geq0$, we see that the initial condition $x(0)=0$ selects the equation (\ref{2.8}) over (\ref{2.9}) for $t$ near $t=0$.
In other words, 
 a solution always starts its evolution following (\ref{2.8}) rather than (\ref{2.9}).

Let $x(t)$ be a local solution of (\ref{2.8}) starting from $x(0)=0$. It will be useful to estimate time $t=t_1>0$ needed for $x(t)$ to cross the level $x=x_1/2$ (say). To this end,
we use
\be\label{2.10}
0\leq x(t)\leq \frac{x_1}2,\quad 0\leq t\leq t_1
\ee
and (\ref{2.8}) to arrive at the inequality
\be\label{2.11}
\dot{x}\geq \sqrt{x}\sqrt{\frac{x_1}{2(\xi+1)}\bigg(x_2-\frac{x_1}2\bigg)},\quad 0\leq t\leq t_1.
\ee
Integrating (\ref{2.11}) over $0\leq t\leq t_1$, we find the following upper bound for $t_1$,
\be\label{2.12}
t_1\leq 2\sqrt{2}\sqrt{\frac{\xi+1}{2x_2-x_1}}.
\ee

We next prove that there is a later time $\ts>t_1$ such that $x(\ts)=x_1$. In \cite{LE}, $\ts$ is called the stagnation time corresponding to the stagnation position $\xs=x_1$
since the velocity of the motion vanishes at $t=\ts$ when the top electrode stretches to the position $x=\xs$.

Indeed, for $t\geq t_1$ so that $x(t)$ stays below $x_1$, that is,
\be\label{2.13}
\frac{x_1}{2}\leq x(t)< x_1 \,\, (=\xs),
\ee
we obtain from (\ref{2.8}) the differential inequality
\be\label{2.14}
\dot{x}\geq \sqrt{\frac{x_1}{2(\xi+1-\frac{x_1}{2})}(x_2-x_1)}\sqrt{x_1-x},\quad t\geq t_1.
\ee
Integrating (\ref{2.14}), we have
\be\label{2.15}
2\bigg(\sqrt{\frac{x_1}2}-\sqrt{x_1-x(t)}\bigg)\geq \sqrt{\frac{x_1}{2(\xi+1-\frac{x_1}{2})}(x_2-x_1)}(t-t_1),\quad t\geq t_1.
\ee
Such a relation establishes the expected property that there must be a finite $\ts>t_1$ such that
\be\label{2.16}
x(\ts)=\xs=x_1
\ee
as claimed.

\section{Periodic solutions}

We consider what happens after the stagnation time $\ts$. Naturally, we may intend to use (\ref{2.8}) again to extend the solution. However, this cannot be done because
(\ref{2.8}) results in an increasing solution but the right-hand side of (\ref{2.8}) is invalid for $x>\xs=x_1$. Therefore, we need the solution to decrease after
$\ts$. Fortunately, such a feature is reflected in (\ref{2.9}). In other words, we are led to use (\ref{2.9}) to extend the solution after $\ts$. In fact, we can do
so explicitly by setting
\be\label{3.1}
x(t)\to x(2\ts-t),\quad t\geq \ts.
\ee
In other words, the definition
\be\label{3.2}
\tilde{x}(t)=\left\{\begin{array}{cc} x(t),& 0\leq t\leq \ts,\\
                                     x(2\ts-t), & \ts\leq t\leq 2\ts\end{array}\right.
\ee
gives us a solution of (\ref{Eq}) which satisfies (\ref{2.8}) for $0\leq t\leq \ts$ and (\ref{2.9}) for $\ts\leq t\leq 2\ts$. This construction clearly produces a
periodic solution with period
\be\label{3.3}
\tp=2\ts.
\ee

In view of (\ref{2.12}) and (\ref{2.15}), we arrive at an upper bound estimate for $\ts$,
\be\label{3.4}
\ts\leq 2\bigg(\sqrt{\frac{\xi+1-\frac{x_1}2}{x_2-x_1}}+\sqrt{2}\sqrt{\frac{\xi+1}{2x_2-x_1}}\bigg).
\ee


\section{Non-periodic or touch-down solution}

We now study what happens when (\ref{2.3}) is violated.

(i) We first assume that
\be\label{4.1}
V^2>\frac14 (\xi+1)^3.
\ee
 Since (\ref{4.1}) implies that the quadratic function $q(x)$ defined in Section 3 is
positive definite,
\be\label{4.2}
\min\{q(x)\} =\frac{V^2}{\xi+1}-\frac{(\xi+1)^2}4\equiv a^2>0.
\ee
Hence (\ref{2.8}) leads to
\be\label{4.3}
\dot{x}=\sqrt{\frac{x}{\xi+1-x} q(x)}\geq \frac{a}{\sqrt{\xi+1}}\sqrt{x},\quad t\geq 0.
\ee
Integrating (\ref{4.3}) over $[0, t]$, we find
\be\label{4.4}
2\sqrt{x(t)}\geq \frac a{\sqrt{\xi+1}} t,
\ee
which indicates that there must be a finite time $\tc>0$, called the ``contact time",  such that
\be\label{4.5}
x(\tc)=\lim_{t\to \tc} x(t)=1.
\ee
Beyond $\tc$, the solution does not exist and (\ref{4.5}) since now the top electrode is in contact with the
dielectric layer. In particular, when $\xi=0$ so that the bottom electrode is not coated with
dielectric, (\ref{4.5}) implies that the switch device is ``short-circuited". Thus,
the critical level
\be\label{4.6}
V^2=\frac14 (\xi+1)^3,
\ee
which borders the dynamic touch-down and periodic response situations of the device, is rightfully referred to as the pull-in voltage \cite{LE}.

(ii) We now assume (\ref{4.6}). From (\ref{2.4}), we see that 
\be\label{4.7}
q(x)=\bigg(x-\frac12(\xi+1)\bigg)^2
\ee
and (\ref{2.8}) becomes
\be\label{4.8}
\dot{x}=\sqrt{\frac x{\xi+1-x}}\,\,\bigg|x-\frac12(\xi+1)\bigg|.
\ee

Since $\xi\geq0$, we see that $\xi+1>(\xi+1)/2$ and the right-hand side of (\ref{4.8}) is Lipschitzian around $x=(\xi+1)/2$. As a consequence,
the uniqueness of the solution to the initial value problem of the ordinary differential equation (\ref{4.8}) holds in a neighborhood of $x=(\xi+1)/2$. Since $x=(\xi+1)/2$ is
an equilibrium point of (\ref{4.8}), it is not attainable at a finite time $t>0$ by a solution different from this equilibrium point.
Consequently, we must have
\be\label{4.9}
\lim_{t\to\infty} x(t)=\frac12(\xi+1),
\ee
which is the limiting position coordinate of the top electrode at the pull-in voltage given in (\ref{4.6}) and referred to as the pull-in position coordinate \cite{LE} of the
top electrode.

\section{A precise dynamic response theorem} 

In summary of the study carried out in the previous sections, following \cite{LE}, a complete description of the dynamic response of an undamped electrostatic MEMS actuator 
subject to a linear elastic force is obtained as follows.

\begin{theorem}\label{theorem1}
Consider the differential equation (\ref{Eq}) governing the dynamics of an electrostatic actuator with a top movable electrode subject to a linear elastic force and the free
initial condition (\ref{IC}). When the applied voltage is low such that $V<\frac12 (\xi+1)^{3/2}$, the response of the top electrode is
periodic characterized by a stagnation time $\tts>0$ at which the electrode stretches to its maximum distance
\be\label{4.10}
\xxs=\frac12(\xi+1)-\frac12\sqrt{(\xi+1)^2-\frac{4V^2}{(\xi+1)}},
\ee
and the electrode follows the same route to return to its original free initial position at $t=2\tts$. In other words, the response is 
represented by a period $\ttp=2\tts$ function which is symmetric about the point $t=\tts$. When the applied voltage is high such
that $V>\frac12 (\xi+1)^{3/2}$, the response of the top electrode is such that, it is driven into contact with the dielectric layer coated over
the bottom electrode plate or the bottom electrode plate itself if it is not coated, at a finite time
$\tc>0$. In other words, $x(t)\to1$ as $t\to \ttc$. As a consequence, the critical pull-in voltage is identified explicitly as
$\VV=\frac12 (\xi+1)^{3/2}$. When the applied voltage is at the level $\VV$, the plate-shaped electrodes will never
make contact with each other but the top electrode approaches an equilibrium position, or the pull-in position, characterized by $x(t)\to
\xxd=\frac12(\xi+1)$ as $t\to\infty$.
\end{theorem}

In the next section, we will see how the method developed in proving Theorem \ref{theorem1} for a linear-elastic force situation can be adapted to obtain
an equally sharp dynamic response theorem for a widely useful nonlinear-elastic force situation.

\section{A nonlinear elastic force situation} 

We now consider the undamped parallel-plate actuator subject to a nonlinear elastic force. A typical model is given
by the Hamiltonian \cite{CO,FM}
\be \label{H3a}
H=\frac12 m{\dot{x}}^2+\frac12 k x^2+\frac14 k_3 x^4-\frac{\varepsilon_0 AV^2}{2\left(\frac{d_0}{\varepsilon_\rr}+g -x\right)},
\ee
where $k_3>0$ is usually a small constant which gives rise to a cubicly nonlinear elastic force. Using the new variables and parameters following (\ref{V}) and setting
\be \label{kappa}
\kappa=\frac{k_3 g^2}k,
\ee
the Hamiltonian (\ref{H3a}) is normalized to take the form
\be \label{H3}
H=\frac12 {\dot{x}}^2+\frac12  x^2+\frac14 \kappa x^4-\frac12\frac{V^2}{\left(\xi+1 -x\right)},
\ee
so that the momentum equation governing the motion of the top electrode reads
\be\label{Eq3}
{\ddot{x}}=-\frac{\pa H}{\pa x}=-x-\kappa x^3+\frac12\frac{V^2}{(\xi+1-x)^2},\quad t>0.
\ee
As before, in view of the energy conservation, we can obtain the first integral of (\ref{Eq3}) as
\be \label{FI2}
\dot{x}^2=\frac{x}{\xi+1-x}\left(\frac{V^2}{\xi+1}-(\xi+1-x)x-\frac\kappa 2 (\xi+1-x)x^3\right).
\ee 
For the normalized coated dielectric layer of thickness $\xi\geq0$, the solution of interest satisfies
\be \label{Interval}
0\leq x\leq1.
\ee
Thus, we need to require
\be 
g(x)=\frac{V^2}{\xi+1}-(\xi+1-x)x-\frac\kappa 2 (\xi+1-x)x^3\geq0
\ee
for consistency in (\ref{FI2}) and (\ref{Interval}). Motivated by the study in the previous sections where $\kappa=0$, we shall maintain the global strict convexity of $g(x)$ by imposing
\be \label{kappa2}
\kappa<\frac{16}{3(\xi+1)^2},
\ee 
which is compatible with the fact that $k_3$ is small (cf. (\ref{kappa})) in applications.

We need to investigate when the global minimum of $g(x)$ becomes negative. We can check that
\be 
g'(0)=-(\xi+1)<0,\quad g'(\xi+1)=(\xi+1)\left(\frac\kappa2(\xi+1)^2+1\right)>0.
\ee
Hence there is a unique point $x_0\in \bfR$ such that
\be \label{g'}
g'(x_0)=2\kappa x_0^3-\frac32\kappa(\xi+1)x_0^2+2x_0-(\xi+1)=0,\quad 0<x_0<\xi+1,
\ee
since $g(x)$ is strictly convex due to (\ref{kappa2}). Such an $x_0$ is the global minimum point of the function $g(x)$. It is interesting to note that $x_0$ depends only on $\xi$ and $\kappa$ but not on
$V$, which enables us to treat
the following three cases individually in a nice way.

Case (i): $V$ is sufficiently low so that $g(x_0)<0$. Since it is clear that
\be 
g(x)>0,\quad x\in (-\infty,0]\cup[\xi+1,\infty),
\ee
we see that there are exactly two points, say $x^*_1$ and $x^*_2$, such that
\be 
g(x^*_1)=g(x^*_2)=0,\quad x^*_1\in(0,x_0),\quad x^*_2\in (x_0,\xi+1).
\ee 
Therefore $g(x)$ enjoys the factorization
\be \label{g}
g(x)=(x_1^*-x)(x^*_2-x)q(x),\quad \forall x,
\ee
where $q(x)$ is a positive-valued quadratic polynomial. Inserting (\ref{g}) into (\ref{FI2}), we arrive at
\be \label{FI3}
\dot{x}^2=\frac{x}{\xi+1-x}(x_1^*-x)(x^*_2-x)q(x).
\ee 
Consequently, we may duplicate the methods in Sections 3 and 4 to show that there is a finite time $\ts^*$ such that a solution satisfying the initial condition (\ref{IC}) climbs to
the stagnation level $\xs^*=x_1^*$ when $t=\ts^*$,
\be 
x(\ts^*)=\xs^*=x_1^*,
\ee
and descends back to the initial position $x=0$ after the same amount of elapse of time, $\ts^*$. In other words, we obtain a symmetric periodic solution of period $\tp^*=2\ts^*$.

It will be interesting to compare $\xs^*$ with $\xs$ obtained in Section 4, when $V<\V$. In fact, with the notation of Section 4, we have
\bea 
g(x_1)&=&-\frac\kappa2(\xi+1-x_1)x_1^3<0,\label{gx1-}\\
g'(x_1)&=&-\frac\kappa2 x_1^2(\xi+1)-(1+\kappa x_1^2)\sqrt{(\xi+1)^2-\frac{4V^2}{\xi+1}}<0.\label{g'x1-}
\eea
From (\ref{g'x1-}) and $g''(x)>0$, we have $x_1<x_0$; from (\ref{gx1-}) and $g'(x)<0$ for $x\in[0,x_0)$, we conclude that
\be 
\xs^*=x^*_1<x_1=\xs=\frac12(\xi+1)-\frac12\sqrt{(\xi+1)^2-\frac{4V^2}{(\xi+1)}},
\ee 
whenever $V<\V$. That is, in the situation of a nonlinear elasticity for the top electrode, the stagnation distance is shorter than that of the linear situation.

Case (ii): $V$ is sufficiently high so that $g(x_0)>0$. We may use the method in Section 5 to show that there is a finite contact time, $\tc^*$, such that
\be 
x(\tc^*)=\lim_{t\to\tc^*}x(t)=1.
\ee
In other words, the top electrode plate stops at time $t=\tc^*$ to take a touch-down contact position with the bottom electrode plate.

Case (iii): $V$ assumes a critical value, say $\Vc$, so that $g(x_0)=0$. That is, $x_0$ satisfies
\be \label{Vc}
\frac{{\Vc}^2}{\xi+1}=(\xi+1-x_0)x_0+\frac\kappa 2 (\xi+1-x_0)x_0^3.
\ee
Since $x_0$ also satisfies (\ref{g'}), it is a double root of the equation $g=0$, which implies the factorization
\be \label{gd}
g(x)=(x_0-x)^2 q(x),
\ee
where $q(x)$ is a positive-valued quadratic polynomial. Substituting (\ref{gd}) into (\ref{FI2}), we have
\be \label{Eqx0}
\dot{x}=\sqrt{\frac {x q(x)}{\xi+1-x}}|x_0-x|.
\ee
Again, we may use the method in Section 5 to show that a solution of (\ref{Eqx0}) satisfying the initial condition (\ref{IC}) is defined for all $t>0$ and has the monotone limit
\be 
\lim_{t\to\infty} x(t)=x_0.
\ee 

Following Section 5, we define the critical voltage $\Vc$ given in (\ref{Vc}) to be the dynamic pull-in voltage, $\V^*$, that is
\be \label{V*}
\V^*=\sqrt{(\xi+1)\left((\xi+1-x_0)x_0+\frac\kappa 2 (\xi+1-x_0)x_0^3\right)},
\ee 
and the corresponding asymptotic position coordinate $x_0$, given in (\ref{g'}), to be the dynamic pull-in position coordinate, $\xd^*$, 
\be 
\xd^*=x_0.
\ee 
 However, since
\be 
g'(\xd)=g'\left(\frac12(1+\xi)\right)=-\frac\kappa8(\xi+1)^3<0,
\ee
and $g''(x)>0$, we see that there holds
\be \label{x0xd}
x_0=\xd^*>\xd,\quad \quad\mbox{if }\kappa>0.
\ee

Furthermore, use $h(x)$ to denote the function
\be \label{h}
h(x)=(\xi+1-x)x+\frac\kappa2(\xi+1-x)x^3,
\ee
so that the right-hand side of (\ref{Vc}) is simply $h(x_0)$. Then $h'(x)=-g'(x)$. In particular, $h'(x)>0$ for $x<x_0$, which implies that $h(\xd^*)>h(\xd)$.
An immediate consequence of this fact, combined with (\ref{Vc}), is that
\be 
\frac{(\V^*)^2}{\xi+1}>\frac{\V^2}{\xi+1}+\frac\kappa{32}(\xi+1)^4,\quad\mbox{if }\kappa>0,
\ee
in particular, 
\be 
\V^*>\V,\quad\mbox{if }\kappa>0,
\ee
which is a naturally expected result, which says that the presence of a cubicly nonlinear elastic force for the motion of the movable electrode of the actuator enhances both dynamic pull-in voltage and
dynamic pull-in position coordinate.

It is worth noting that, in the limiting situation when $\kappa=0$, the formulas (\ref{g'}) and (\ref{V*}) recover the linear elasticity results obtained earlier, so that
\be 
\xd^*\to\xd=\frac12(\xi+1),\quad \V^*\to\V=\frac12(\xi+1)^{\frac32},\quad \kappa\to0.
\ee

The study just carried out enables us to arrive at the following second sharp dynamic response theorem concerning the MEMS actuator subject to a cubicly nonlinear elastic force.

\begin{theorem}
Consider the differential equation (\ref{Eq3}) describing the dynamics of an undamped electrostatic actuator subject to a cubicly nonlinear elastic force and the free
initial condition (\ref{IC}), governed by the normalized Hamiltonian (\ref{H3}), and assume that the nonlinear elasticity is weak so that
(\ref{kappa2}) holds. Then the full picture of the dynamic response of the system is completely determined by the unique solution 
$ 
\xd^*=x_0
$
of the cubic equation (\ref{g'}), which necessarily lies in the interval $(0,\xi+1)$ and defines the dynamic pull-in voltage $\V^*$ through the formula (\ref{V*}).
When the applied voltage is low such that $V<\V^*$, the response of the actuator is
periodic characterized by a stagnation time $\ts^*>0$ at which the top movable electrode stretches to its maximum distance $\xs^*$ so that
\be\label{}
\xs^*<\xs=\frac12(\xi+1)-\frac12\sqrt{(\xi+1)^2-\frac{4V^2}{(\xi+1)}},\quad\mbox{if }V<\V,
\ee
and the top electrode follows the same route to return to its initial position at $t=2\ts^*$. In other words, the response is 
represented by a period $\tp^*=2\ts^*$ function which is symmetric about the point $t=\ts^*$. 
Furthermore, as a function of $\kappa$ and $V$, the stagnation position coordinate $\xs^*$ is smooth, strictly decreasing with respect to $\kappa$, and increasing
with respect to $V$.
When the applied voltage is high such
that $V>\V^*$, the response of the actuator is so that, the applied voltage drives the top electrode plate into contact with the dielectric layer coated on the
bottom electrode plate or the bottom electrode plate when there is no coated dielectric layer, at a finite time
$\tc>0$. In other words, $x(t)\to 1$ as $t\to \tc$.  When the applied voltage is at the level $\V^*$, the top electrode plate will never
make contact with the dielectric layer or the bottom electrode plate but approach an equilibrium position, or the pull-in position, characterized by $x(t)\to
\xd$ as $t\to\infty$. Moreover, $\xd^*$ and $\V^*$ depend on $\kappa$ smoothly and are strictly increasing, so that they converge to $\xd$ and $\V$ in the zero nonlinear-elastic force limit
$\kappa\to0$.
\end{theorem}

\begin{proof} We only need to prove the smoothness and monotonicity of $\xs^*$, $\xd^*$, and $\V^*$ with respect to the designated variables.

To proceed, assume $V<\V^*$. In this situation, we first notice that the condition (\ref{kappa2}) ensures that $g'(\xs^*)=g'(x^*_1)<0$ since $g'(x_0)=0$ and $g''(x)>0$ for all $x$.
Thus, the implicit equation
\be 
G(x_1^*,\kappa,V)\equiv\frac{V^2}{\xi+1}-(\xi+1-x_1^*)x_1^*-\frac\kappa 2 (\xi+1-x_1^*)(x_1^*)^3=g(x_1^*)=0
\ee
gives us
\be \label{GP}
\frac{\pa G}{\pa x_1^*}(x_1^*,\kappa,V)=g'(x_1^*)<0,
\ee
and the smooth dependence of $x_1^*$ on $\kappa$ and $V$ follows as a result of the implicit function theorem. Besides, after implicit differentiation, we
have 
\bea 
\left(\frac{\pa G}{\pa x_1^*}\right)\frac{\pa x_1^*}{\pa \kappa}&=&\frac12(\xi+1-x_1^*)(x_1^*)^3>0,\label{G1}\\
\left(\frac{\pa G}{\pa x_1^*}\right)\frac{\pa x_1^*}{\pa V}&=&-\frac{2V}{\xi+1}<0.\label{G2}
\eea
Therefore, combining (\ref{GP})--(\ref{G2}), we arrive at the conclusion
\be 
\frac{\pa x_1^*}{\pa \kappa}<0,\quad \frac{\pa x_1^*}{\pa V}>0.
\ee

We now show that $\xd^*$ and $\V^*$ are smoothly defined increasing functions of the parameter $\kappa$.
To this end, recall that $\xd^*=x_0$ satisfies the equation (\ref{g'}), which may be rewritten as $F(x_0,\kappa)=0$. Of course,
\be \label{Fx0}
\frac{\pa F}{\pa x_0}=6\kappa\left(x_0^2-\frac12(\xi+1)x_0+\frac1{3\kappa}\right)>0,
\ee
due to the condition (\ref{kappa2}), which establishes that $x_0$ is a smooth function of $\kappa$ determined implicitly by $F(x_0,\kappa)=0$. Using (\ref{g'}), we have
\bea \label{PF}
\left(\frac{\pa F}{\pa x_0}\right)\frac{\dd x_0}{\dd\kappa}&=&-\left(\frac{\pa F}{\pa\kappa}\right)=-\left(2x_0^3-\frac32(\xi+1)x_0^2\right)\nn\\
&=&\frac2\kappa\left(x_0-\frac12(\xi+1)\right)>0,
\eea
in view of (\ref{x0xd}) and the fact that $\xd=\frac12(\xi+1)$ obtained earlier. Combining (\ref{Fx0}) and (\ref{PF}), we have $\frac{{\small{\dd} }x_0}{{\small{\dd}}\kappa}>0$ as claimed.

With the notation (\ref{h}), we rewrite (\ref{V*}) as $\V^*=\sqrt{(\xi+1)h(x_0)}$. Therefore
\be 
\frac{\dd\V^*}{\dd \kappa}=\frac12\sqrt{\frac{(\xi+1)}{h(x_0)}}\left(h'(x_0)\frac{\dd x_0}{\dd\kappa}+\frac12(\xi+1 -x_0)x_0^3\right)>0,
\ee
since $h'(x_0)=-g'(x_0)=0$ and $x_0<\xi+1$.

The proof of the theorem is complete.
\end{proof}

The physical meaning of the monotone dependence of $\xs^*=x_1^*$ on $V$ and $\kappa$ is natural to expect since a larger voltage $V$ drives the top electrode plate down further but stronger elasticity characterized by
a larger $\kappa$ parameter gives greater stiffness, which in turn leads to a smaller stretching distance for the periodic motion.

We also note that our method may be used to study the model governed by the following normalized extended Hamiltonian
\be \label{H4}
H=\frac12 {\dot{x}}^2+\Phi(x)-\frac12\frac{V^2}{\left(\xi+1 -x\right)},
\ee
where $\Phi(x)\geq0$ with $\Phi(0)=0$ is a general elastic potential density function and damping is absent. In fact, since the energy conservation law gives us
the relation
\be 
\dot{x}^2=\frac{V^2x}{(\xi+1)(\xi+1-x)}-2\Phi(x),
\ee
we see that the crucial condition to be imposed so that our method works would be
the strict convexity of the function
\be 
\Psi(x)=\frac{V^2}{\xi+1}-2(\xi+1-x)\frac{\Phi(x)}x,\quad x>0.
\ee 


\section{Touch-down as a universal dynamical property} 

In this section, we shall identify some simple conditions  under which
the touch-down phenomenon inevitably takes place when a coupling parameter resembling the applied voltage is sufficiently large. For greater generality, we will include damping in the model
so that the normalized equation of motion assumes the form
\be \label{Eq8}
\ddot{x}+\mu\dot{x}+f(x,t)=\lm g(x,t),
\ee 
where $\mu\geq 0$ is a damping constant, a physically meaningful solution $x=x(t)$ stays in the interval $[0,a]$ for some $a>0$, and $f$ and $g$ are functions of $x$ and $t$ so that
there are positive constants $C_1$ and $ C_2$ such that the bounds
\bea 
\sup\{|f(x,t)|\,|\,x\in[0,a], t\geq0\}&\leq& C_1,\label{C1}\\
\inf \{g(x,t)\,|\, x\in[0,a), t\geq0\}&\geq& C_2.\label{C2}
\eea
hold.

We will assume $a$ to be the touch-down position for the solution so that the function $g$, a term resembling the Coulomb force, may be singular at $x=a$.
We have

\begin{theorem} If $\lm>0$ is sufficiently large, there is an associated finite time $\tc>0$, depending on $\lm$ and $a$, such that the solution of
(\ref{Eq8}) subject to the initial condition (\ref{IC}) monotonically climbs to the touch-down position at $t=\tc$. That is,
\be \label{xtc}
x(\tc)=\lim_{t\to\tc}x(t)=a\quad \mbox{and}\quad \dot{x}(t)>0\quad\forall t\in(0,\tc).
\ee
\end{theorem}
\begin{proof} First assume $\mu>0$. In view of (\ref{Eq8})--(\ref{C2}) and the initial condition (\ref{IC}), we have, after multiplying (\ref{Eq8}) by $\e^{\mu t}$ and integrating and using the
condition $0\leq x(t)<a$,
the relation
\bea 
\e^{\mu t}\dot{x}(t)&=&\int_0^t\e^{\mu \tau}\left(\lm g(x(\tau),\tau)-f(x(\tau),\tau)\right)\,\dd\tau\nn\\
&\geq&\frac1\mu(\lm C_2-C_1)(\e^{\mu t}-1),
\eea
which establishes $\dot{x}(t)>0$ for all $t>0$ provided that $\lm C_2>C_1$ and leads to the lower estimate
\be 
x(t)\geq\frac1\mu(\lm C_2-C_1)\left(t-\frac1\mu(1-\e^{-\mu t})\right),\quad t\geq0.
\ee
Hence there is a finite time $\tc>0$ such that (\ref{xtc}) holds when $\lm C_2>C_1$.

Next, if $\mu=0$, then an integration of (\ref{Eq8})  simply gives us $\dot{x}(t)\geq (\lm C_2-C_1)t$ as before and the rest follows obviously.

 Thus the theorem is proven.
\end{proof}

\medskip 

In summary,  we have presented two theorems which give precise descriptions of the global dynamic behavior of an electrostatic undamped actuator, subject to linear and cubicly nonlinear
elastic forces, arising in MEMS technology. The method developed is applicable to the study of the models with more general elastic forces.

\end{document}